%% file: WildObjects.tex
\documentclass[
11pt, 
a4paper, 
oneside, 
headinclude,footinclude, 
BCOR5mm, 
]{scrartcl}

\input{structure.tex} 

\title{\normalfont\spacedallcaps{A Survey on wild mathematics}} 

\author{\spacedlowsmallcaps{Raffaella Mulas}\thanks{Max Planck Institute for Mathematics in the Sciences, Leipzig, Germany\newline Email address: raffaella.mulas@mis.mpg.de}} 

\date{} 

\usepackage[margin=3cm]{geometry}
\usepackage{amsthm} 
\usepackage{amssymb}
\usepackage{amsopn}
\usepackage{color}
\usepackage[english]{babel} 
\usepackage{hyperref}
\usepackage{youngtab}
\usepackage{amsmath, amsfonts, amssymb} 
\usepackage{graphicx}
\allowdisplaybreaks 
\usepackage[utf8]{inputenc} 
\usepackage{xfrac}    
\usepackage[draft]{todonotes}
\usepackage{amssymb}

\usepackage{mathtools} 
\mathtoolsset{showonlyrefs,showmanualtags}
\numberwithin{equation}{section}
\usepackage{amsmath}
\usepackage{float}

\theoremstyle{plain}
\newtheorem{theorem}{Theorem}

\theoremstyle{definition}

\newtheorem{definition}[theorem]{Definition}

\theoremstyle{remark}

\newcommand{\Q}{\mathbb{Q}}

\begin{document}


\renewcommand{\sectionmark}[1]{\markright{\spacedlowsmallcaps{#1}}} 
\lehead{\mbox{\llap{\small\thepage\kern1em\color{halfgray} \vline}\color{halfgray}\hspace{0.5em}\rightmark\hfil}} 

\pagestyle{scrheadings} 

\maketitle 




\begin{abstract}
    We usually construct mathematical objects that are accessible, on which we can put our hands, but a huge part of the mathematical existing is actually wild. Here we explore part of the wild world: its inhabitants are knots that are infinitely knotted, spheres that try to hug themselves, colorful earrings, labyrinths of labyrinths...
\end{abstract}
\begin{figure}[H]
		\begin{center}
			\includegraphics[width=7cm]{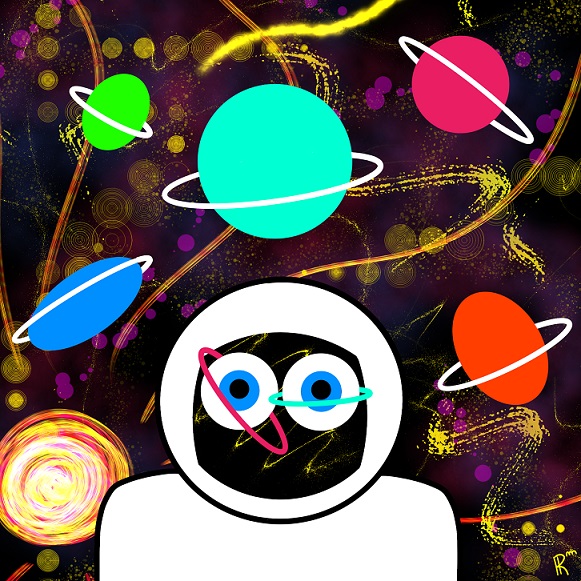}
		\end{center}
		\caption{Sense of Monstrosity by RM.}\label{Figuremonst}
	\end{figure}

Let us give a general definition of \emph{wild} based on the human nature, before exploring some wild objects. During this exploration we shall give more specific and formal definitions of \emph{wild}.\newline
\begin{definition}\label{defwild}
We say that something is \emph{wild} if it gives us a sense of \emph{Monstrosity}, i.e. if we have the feeling that it is inaccessible and we therefore feel puzzled and overwhelmed in front of it (Figure \ref{Figuremonst}). We say that something is \emph{tame} if it's not wild.
\end{definition}

	\newpage
	\section{Irrational numbers}
The elements of $\Q$ are tame. In fact, we can conceive the rational numbers and we can even write them with a computer. On the other hand, the irrational numbers are wild because, in order to describe them, we need an infinite number of decimal digits that are not periodic.

\section{Wild knots}
\begin{figure}[H]
		\begin{center}
			\includegraphics[width=13cm]{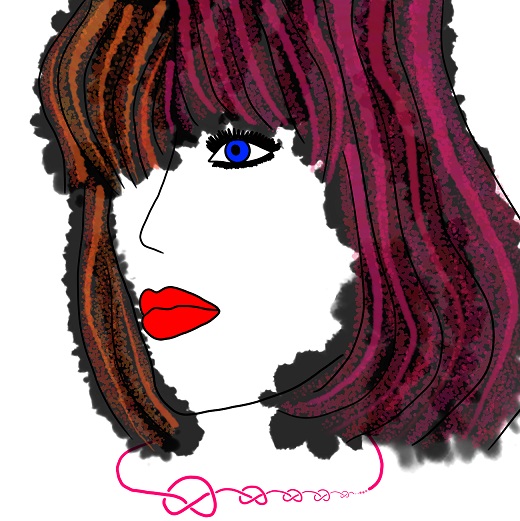}
		\end{center}
		\caption{Wild knot necklace by RM.}\label{FigureKnot}
	\end{figure}

We define a \emph{knot} as a topological embedding of $S^1$ into $\mathbb{R}^{3}$. We say that a knot is \emph{tame} if it is isotopic to a finite closed polygonal chain. We say that a knot is \emph{wild} if it's not tame (Figure \ref{FigureKnot}). Note at this point that the definition of tame and wild knot is in line with Definition \ref{defwild}.\newline

\newpage
\section{Alexander sphere}
Related to wild knots, there is an interesting pathological example of wildness discovered by  J. W. Alexander: it is called the \emph{Alexander horned sphere}, see \cite{Alexander, Hatcher}. This wild topological space is a subspace $A\subset\mathbb{R}^3$ homeomorphic to $S^2$ such that the unbounded component of $\mathbb{R}^3\setminus A$ is not simply-connected as it is for the standard $S^2\subset \mathbb{R}^3$.
	\begin{minipage}{7cm}
It can be constructed as follows:
\begin{itemize}
    \item Consider a standard torus and remove a radial slice from it.
    \item To each side of the cut, connect a standard torus,  interlinked with the torus on the other side.
    \item Repeat the previous steps infinitely many times.
\end{itemize}The limit object obtained with this construction is the Alexander sphere. 
\end{minipage}
\begin{minipage}{7cm}
\begin{figure}[H]
		\begin{center}
			\includegraphics[width=7cm]{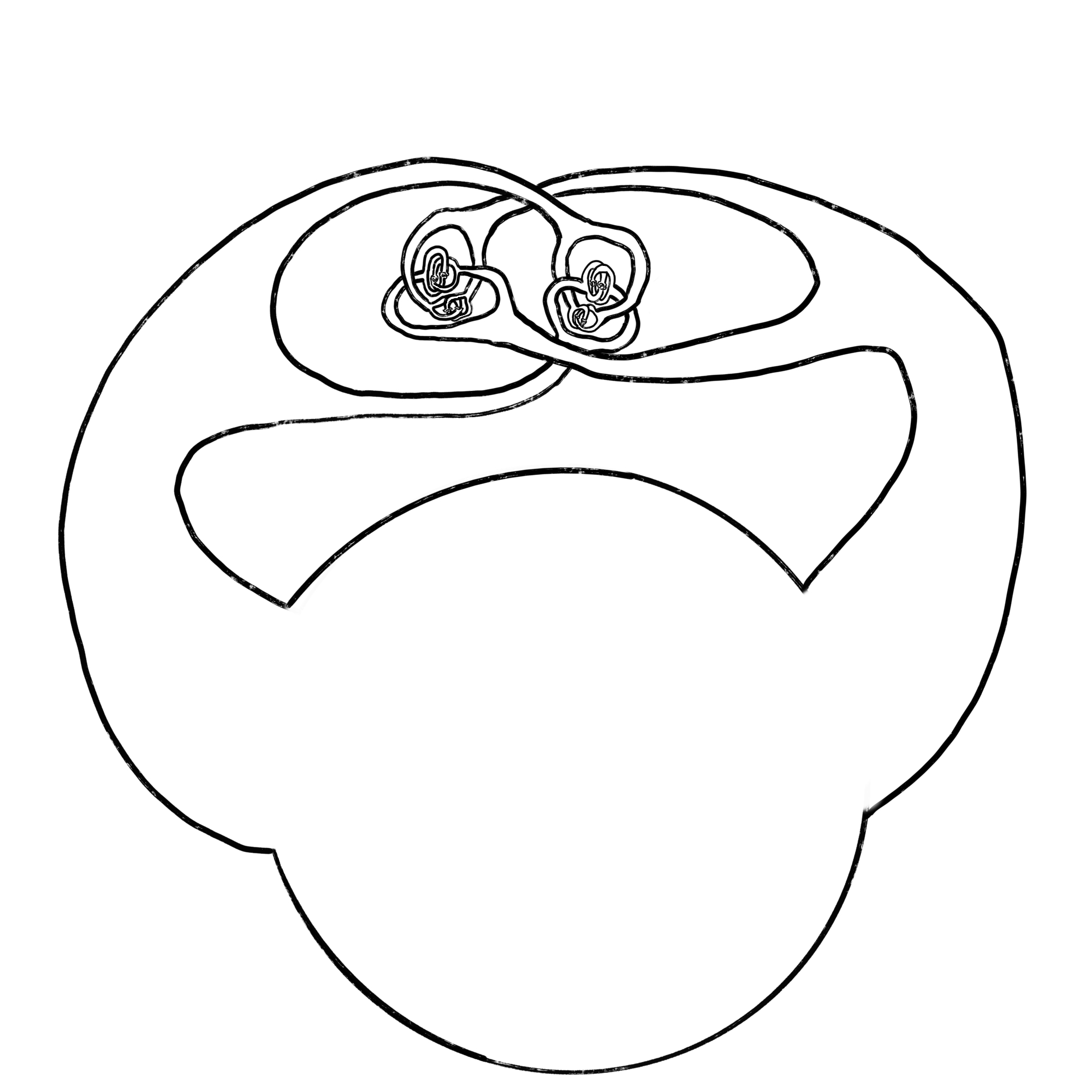}
		\end{center}
		\caption{Wild Alexander sphere by RM.}\label{FigureAlex1}
	\end{figure}
\end{minipage}\newline
\begin{minipage}{7cm}
\begin{figure}[H]
		\begin{center}
			\includegraphics[width=7cm]{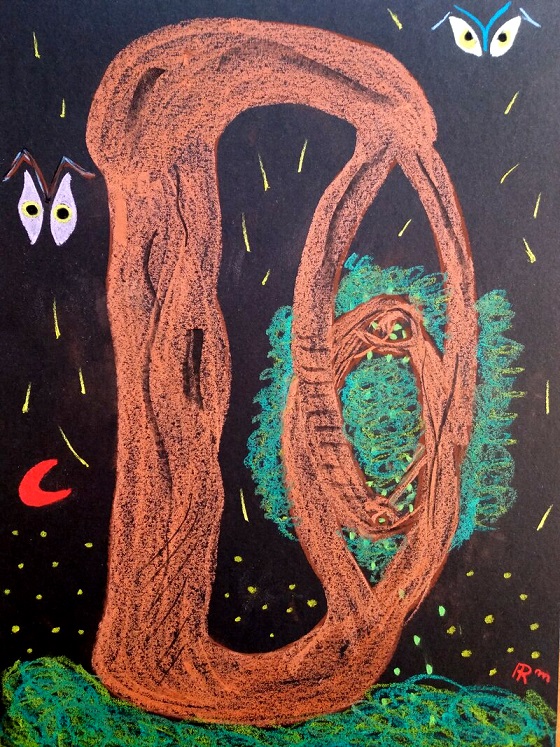}
		\end{center}
		\caption{Very wild Alexander sphere by RM.}\label{FigureAlex2}
	\end{figure}
\end{minipage}\hfill
	\begin{minipage}{7cm}
 It is a sphere that tries to hug itself and it is wild in the sense of Definition \ref{defwild}. Here we propose two representations of the Alexander sphere: Figure \ref{FigureAlex1} clearly shows its inductive construction; Figure \ref{FigureAlex2} shows its wild nature.
	\end{minipage}
\section{The Cantor set}
\begin{minipage}{7cm}
	The Cantor set (Figure \ref{FigureCantor}) is constructed iteratively as follows:
\begin{itemize}
    \item From the unit interval $[0,1]$, remove the open central third $(\frac{1}{2},\frac{2}{3})$.
    \item From each remained interval, remove the open central third.
    \item Repeat this process infinitely many times.
\end{itemize}The Cantor set is the limit object of this process and it is proved to be uncountable \cite[Proposition 4.36]{Cantor}. Also, every compact metric space is the continuous image of the Cantor set: this is really too much of Monstrosity!
\end{minipage}
\hfill
\begin{minipage}{7cm}
	\begin{figure}[H]
		\begin{center}
			\includegraphics[width=7cm]{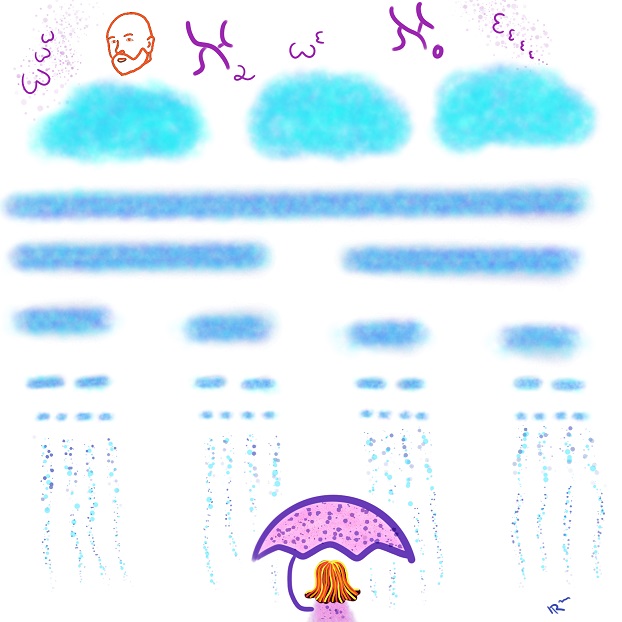}
		\end{center}
		\caption{Cantor rain by RM.}\label{FigureCantor}
	\end{figure}
\end{minipage}

	\section{Hawaiian spheres}\label{ExHawaii}
	\begin{minipage}{8cm}
		\begin{figure}[H]
		\begin{center}
			\includegraphics[width=8cm]{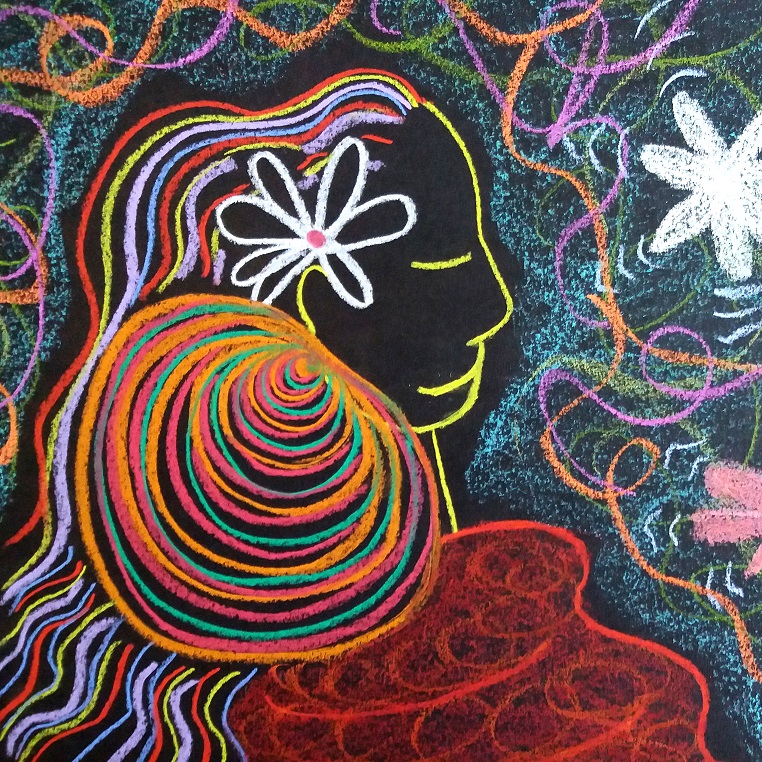}
		\end{center}
		\caption{Hawaiian earrings by RM.}\label{FigureEarrings}
	\end{figure}
	\end{minipage}
\hfill
\begin{minipage}{6cm}
For $m\geq1$, the $m$--dimensional \emph{Hawaiian spheres} are a countable number of $m$--spheres with a single point in common and a metric topology such that the diameter of the spheres tends to zero with increasing index \cite{Milnor}. In Flatland \cite{Flatland}, i.e. in the $2$--dimensional world, the $1$--Hawaiian spheres are also called \emph{Hawaiian earrings} (Figure \ref{FigureEarrings}). The Hawaiian earrings are wild as they are defined by a limit, but the Hawaiian spheres for $m\geq 2$ are even \emph{wilder}. In fact, in this case, these spheres form a $m$--dimensional object with the property that for each $n$ there exists $l>n$ such that the $l$--th homology group is non--zero, as shown in \cite[Theorem 1]{Milnor}.
	\end{minipage}\newline
	
\newpage
\section{Peano curve}
\begin{figure}[H]
		\begin{center}
			\includegraphics[width=11cm]{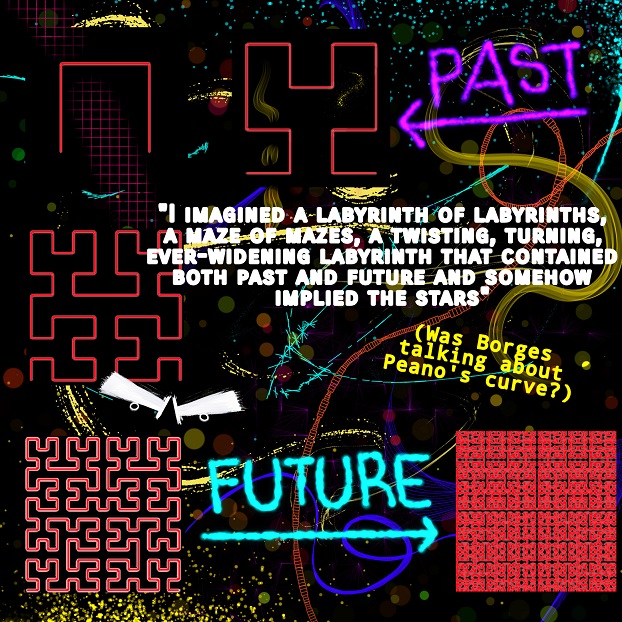}
		\end{center}
		\caption{Peano and Borges by RM.}\label{FigurePeano2}
	\end{figure}

\begin{minipage}{7cm}
A \emph{Peano curve} is a continuous function from the unit interval into the unit square that is surjective (but not injective!). A Peano curve is usually defined as the limit of a sequence of curves, as for example the one constructed by David Hilbert as the limit of the construction in Figure \ref{FigurePeano1}. One can prove that such limit curve exists, is continuous and is surjective. \emph{Wooow!} Not only it gives us a feeling of Monstrosity but, as shown in Figure \ref{FigurePeano2}, it makes us think about \emph{The Garden of Forking Paths} from Borges \cite{Borges}.
\end{minipage}
\hfill
\begin{minipage}{7cm}
\begin{figure}[H]
		\begin{center}
			\includegraphics[width=7cm]{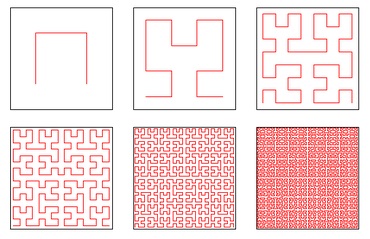}
		\end{center}
		\caption{Six iterations of the Peano curve constructed by David Hilbert (picture from Wikipedia, licensed under the Creative Commons Attribution-Share Alike 3.0 Unported license).}\label{FigurePeano1}
	\end{figure}\end{minipage}

	\newpage
	\section{Fractals}
\begin{figure}[H]	\begin{center}\includegraphics[width=15cm]{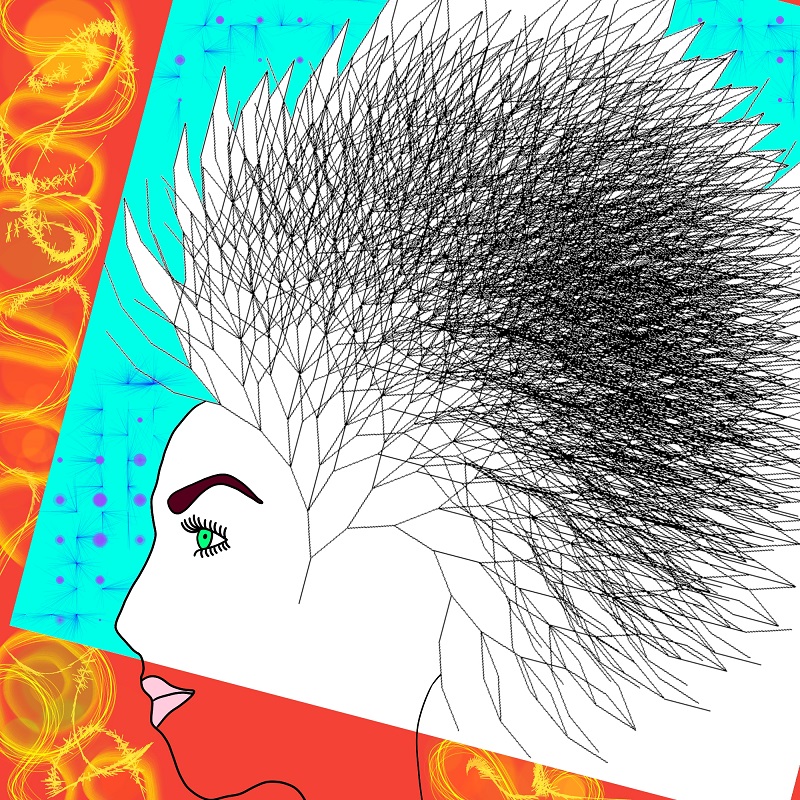}\end{center}\caption{Fractal hair by RM.}\label{FigureHair}\end{figure}

The Alexander sphere, the Cantor set and the Hawaiian spheres are all \emph{fractals}, i.e. objects for which zooming in always reveals the same pattern. Fractals are therefore characterized by self-similarity on different scales and, by Definition \ref{defwild}, they are wild. Fractals also exist in nature: they can be found for instance in coastlines, mountains, trees, clouds, veins and even in the hair, as shown in Figure \ref{FigureHair}.

\newpage

\section{Conclusions}
Wild objects are everywhere around us, in mathematics as well as in other fields, especially if there is some hidden mathematical structure. We have mentioned \emph{The Garden of Forking Paths} of Borges, for instance. This short story appears in Everett's interpretation of quantum mechanics \cite{Quantum}, and the entire literary work of Borges is so much soaked in mathematics that the Argentinian mathematician and writer Guillermo Martínez wrote a book titled \emph{Borges y la Matemática} \cite{Martinez}. We also like to think that Borges' short story about \emph{The Aleph} \cite{Borges2}, described as \emph{one of the points in space that contains all other points}, is a tribute to Cantor's Aleph numbers that measure different \emph{infinities}. How wild is this? \newline

Another writer whose work is dense of mathematics and wildness is Lewis Carroll. We could say for instance that, while the arithmetic operations $+$, $-$, $\cdot$ and $/$ are tame, the operations of Ambition, Distraction, Uglification, and Derision described by Lewis Carroll in \emph{Alice's Adventures in Wonderland} \cite{Alice} are wild. Note that, in the Italian version of the book, the operation of Uglification is called \emph{Mostrificazione}, that can be literally translated with \emph{Monsterification}. This definitely has to do with the sense of Monstrosity! \newline
Also Escher, with his drawings, created wild worlds playing with geometry and, particularly, with  the fact that his objects drawn in two dimensions cannot leave in three dimensions. \newline

\begin{figure}[H]
		\begin{center}
			\includegraphics[width=13cm]{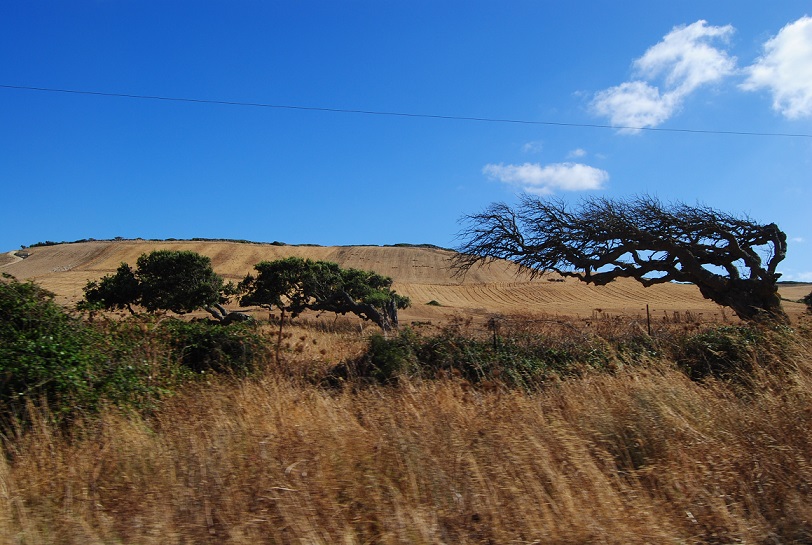}
		\end{center}
		\caption{Wild trees bent over by the wild wind in Sardinia. Picture taken by RM.}\label{FigureWind}
	\end{figure}

Then there is Nature, with its wildness (Figure \ref{FigureWind}). It is not clear, however, whether mathematics has a wildness that reminds of nature, or whether the wildness in nature reminds of wild mathematics.

\subsection*{Acknowledgements}The author would like to thank Antonio Lerario for the intersting discussions that inspired this work, and she is grateful to Eleonora Andreotti and Jürgen Jost for the support and the valuable comments.
	\bibliographystyle{unsrt}
	\bibliography{WildObjects}	

\end{document}

%% file: structure.tex
\usepackage[
nochapters, 
beramono, 
eulermath,
pdfspacing, 
dottedtoc 
]{classicthesis} 

\usepackage{arsclassica} 

\usepackage[T1]{fontenc} 

\usepackage[utf8]{inputenc} 

\usepackage{graphicx} 
\graphicspath{{Figures/}} 

\usepackage{enumitem} 

\usepackage{lipsum} 

\usepackage{subfig} 

\usepackage{amsmath,amssymb,amsthm} 

\usepackage{varioref} 

\hypersetup{
colorlinks=true, breaklinks=true, bookmarks=true,bookmarksnumbered,
urlcolor=webbrown, linkcolor=RoyalBlue, citecolor=webgreen, 
pdftitle={}, 
pdfauthor={\textcopyright}, 
pdfsubject={}, 
pdfkeywords={}, 
pdfcreator={pdfLaTeX}, 
pdfproducer={LaTeX with hyperref and ClassicThesis} 
}

%% file: WildObjects.bbl
\begin{thebibliography}{10}

\bibitem{Alexander}
J.~W. Alexander.
\newblock An example of a simply connected surface bounding a region which is
  not simply connected.
\newblock {\em Proceedings of the National Academy of Sciences of the United
  States of America}, 10(1):8--10, 1924.

\bibitem{Hatcher}
Allen Hatcher.
\newblock {\em Algebraic Topology}.
\newblock Cambridge University Press, 2002.

\bibitem{Cantor}
Steven~G. Krantz.
\newblock {\em Foundations of Analysis}.
\newblock CRC Press, 2014.

\bibitem{Milnor}
M.~G. Barratt and John Milnor.
\newblock An example of anomalous singular homology.
\newblock {\em Proc. Amer. Math. Soc.}, 13:293--297, 1962.

\bibitem{Flatland}
Edwin~A. Abbott.
\newblock {\em Flatland: A Romance of Many Dimensions}.
\newblock Seeley \& Co., 1884.

\bibitem{Borges}
Jorge~Luis Borges.
\newblock {\em The Garden of Forking Paths}.
\newblock Editorial Sur, 1941.

\bibitem{Quantum}
E.~Hugh, J.~A. Wheeler, B.~S. DeWitt, L.~N. Cooper, D.~Van~Vechten, and
  N.~Graham.
\newblock {\em The Many-Worlds Interpretation of Quantum Mechanics}.
\newblock Princeton Series in Physics. Princeton, NJ: Princeton University
  Press, 1973.

\bibitem{Martinez}
Guillermo Martínez.
\newblock {\em Borges y la Matemática}.
\newblock Universidad de Buenos Aires, 2003.

\bibitem{Borges2}
Jorge~Luis Borges.
\newblock {\em The Aleph}.
\newblock Editorial Losada, 1945.

\bibitem{Alice}
Lewis Carroll.
\newblock {\em Alice's Adventures in Wonderland}.
\newblock Macmillan, 1865.

\end{thebibliography}
